\documentstyle{article}
\begin{document}

\centerline{\bf CLOSED AND $\sigma$-FINITE MEASURES}
\centerline{\bf ON THE ORTHOGONAL PROJECTIONS}

\bigskip
\centerline{{MARJAN MATVEJCHUK}
\footnote{Novorossiisk University,
Novorossiisk, Geroev Desantnikov Str., 87, 353922 Russia\,\,
E-mail: Marjan.Matvejchuk@ksu.ru}}

 \bigskip
\hrule
 We characterize the connection between closed and $\sigma$-finite
measures on orthogonal projections of von Neumann algebras.
\hrule

\vskip0.5cm

Let $\cal A$ be a von Neumann algebra acting in a separable complex Hilbert
space $H$ and let ${\cal A}^{pr}$ be the set of all orthogonal projections
(=idempotents) from $\cal A$.
A subset ${\cal M}\subseteq {\cal A}^{pr}$ is said to be {\it ideal of
projections} if:

a) $p\leq q$ where $p\in {\cal A}^{pr}$, $q\in \cal M$ $\Rightarrow$ $p\in
{\cal M}$;
b) $p$, $q\in \cal M$ and $\Vert pq\Vert <1$ $\Rightarrow$ $p\vee q\in \cal
M$;
c) $\sup\{p: p\in {\cal M}\}=I$.


Put ${\cal M}_{p}:=\{q: q\in {\cal M}, q\leq p\}$,
$\forall p\in {\cal A}^{pr}$.
Note that ${\cal A}^{pr}$ is the ideal of projections,
 $0\in {\cal M}_{p}$, $\forall p$, and the
 conditions $1)$, $2)$ are fulfilled on ${\cal M}_{p}$.

A function $\mu :{\cal M}\to [0,+\infty]$
is said to be a measure if $\mu (e)=\sum \mu (e_{i})$ for any representation
 $e=\sum e_{i}$.
Let $\mu_{1}:{\cal M}_{1}\to [0,+\infty ]$
and  $\mu_{2}:{\cal M}_{2}\to [0,+\infty ]$ be measures.
The measure  $\mu_{2}$ is said to be the {\it continuation} of $\mu_{1}$ if
${\cal M}_{1}\subset {\cal M}_{2}$ and $\mu_{1}(p)=\mu_{2}(p)$,
$\forall p\in {\cal M}_{1}$.
A projection $p\in {\cal A}^{pr}$ is said to be: {\it projection of finite
 $\mu$-measure} if $\sup\{q\in {\cal M}_{p}\}=p$  and
$\sup\{\mu (q): q\in {\cal M}_{p}\}<+\infty$;
{\it hereditary projection of finite $\mu$-measure}
if $q$ is the projection of finite $\mu$-
measure for any $q\in {\cal A}^{pr}$, $q\leq p$.

The measure $\mu$ is said to be: {\it finite} if
$\mu (p)<\infty$, $\forall p$; {\it infinite} if there exists $p\in\cal M$
such that
$\mu (p)=+\infty$; {\it fully finite} if $\sup\{\mu (p):p\in {\cal M}\}<+
\infty$; {\it closed} if $\mu$ is finite and $p\in\cal M$ if $p$ is the
hereditary projection of finite $\mu$-measure;
$\sigma$-{\it finite} if
${\cal M}={\cal A}^{pr}$ and
there exists a sequence
$\{p_{n}\}\subset {\cal A}^{pr}$ such that $p_{n}\nearrow I$ and
$\mu (p_{n})<+\infty$, $\forall n$.

The following Proposition will be needed in Theorem 3.

{\bf Proposition 1.} {\it Let $\cal A$ be a finite von Neumann algebra
acting in the separable
Hilbert space $H$ and let ${\cal M}\subseteq {\cal A}^{pr}$ be the ideal of
projections.
Then there exists a sequence
$\{e_{n}\} \subset \cal M$ such that $e_{n}\nearrow I$.}

{\it Proof.} Let $\tau$ be a faithful normal finite trace on ${\cal A}^{+}$.
The proof consists of several steps.

i). Fix $\epsilon >0$. Let us prove that there exists
$p_{\epsilon}\in {\cal M}$ such that $\tau (p_{\epsilon})>\tau (I)-\epsilon$.
By b) and by separability of $H$, there exists a sequence
$\{q_{n}\}_{1}^{\infty}\subset \cal M$ with $\sup\{q_{n}\}=I$.

1). Put $f_{1}:=q_{1}$.

2). Let $q_{2}f_{1}q_{2}=\int_{0}^{1+}\lambda de_{\lambda}^{(2)}$
be the spectral decomposition of $q_{2}f_{1}q_{2}$. Here $e^{(2)}_{\lambda}$
is the left continuous decomposition of identity. Fix $k\in N$.
By the left continuous of $e_{\lambda}^{(2)}$,
there exists $\beta\in (0,1)$ such that $\tau (e_{1}^{(2)}-
e_{\beta}^{(2)})< ({\frac 1 2})^{k+1}$. Put $q^{\prime}_{2}:=
q_{2}\wedge e^{(2)}_{\beta}$. By 1),
$q^{\prime}_{2}\in\cal M$.
By the construction of $\beta$,
$\Vert q_{2}^{\prime}f_{1}\Vert \leq\beta$ $(<1)$. Hence
$f_{1}\vee q_{2}^{\prime}\leq q_{1}\vee q_{2}$ and by b),
$f_{2}:=f_{1}\vee q_{2}^{\prime}\in\cal M$ and $\tau (q_{1}\vee q_{2}-f_{2})
<({\frac 1 2})^{k+1}$.

3). Let $q_{3}f_{2}q_{3}=\int_{0}^{1+}\lambda de_{\lambda}^{(3)}$ be the
spectral decomposition of $q_{3}f_{2}q_{3}$. Let us choose
$\beta\in (0,1)$ such that $\tau (e_{1}^{(3)}-e^{(3)}_{\beta})
<({\frac 1 2})^{k+2}$, again. Put $q_{3}^{\prime}:=q_{3}\wedge
e^{(3)}_{\beta}$. By a), $q^{\prime}_{3}\in\cal M$. By the construction of
$\beta$, we have
$\Vert q^{\prime}_{3}f_{2}\Vert\leq \beta$ $(<1)$. Again by b),
 $f_{3}:=f_{2}\vee q_{3}^{\prime}\in\cal M$. Thus
$$
f_{3}\leq f_{2}\vee q_{3}\leq f_{1}\vee q_{2}\vee q_{3}=q_{1}\vee q_{2}\vee
q_{3}
$$
and
$$
\tau (f_{2}\vee q_{3}-f_{2}\vee q^{\prime}_{3})\leq ({\frac 1 2})^{k+2},
\qquad
\tau (q_{1}\vee g_{2}\vee q_{3}-f_{2}\vee g_{3})< ({\frac 1 2})^{k+1}.
$$
Therefore $\tau (q_{1}\vee q_{2}\vee q_{3}-f_{3})\leq ({\frac 1 2})^{k+1}+
({\frac 1 2})^{k+2}$.

Let us continue the process of construction of $\{f_{n}\}$
by the induction with respect to $n$.

n). Let the projection $f_{n-1}\in\cal M$ it was chosen.
 Let $q_{n}f_{n-1}q_{n}=\int_{0}^{1+}\lambda de_{\lambda}^{(n)}$ be
the spectral decomposition of $q_{n}f_{n-1}q_{n}$. Let us
choose $\beta\in (0,1)$
such that $\tau (e_{1}^{(n)}-e^{(n)}_{\beta})
<({\frac 1 2})^{k+n-1}$. Put $q_{n}^{\prime}:=q_{n}\wedge e^{(n)}_{\beta}$.
By a), $q^{\prime}_{n}\in\cal M$. By the construction of
$\beta$, we have
$\Vert q^{\prime}_{n}f_{n-1}\Vert\leq \beta$ $(<1)$. By b),
$f_{n}:=f_{n-1}\vee q_{n}^{\prime}\in\cal M$. Thus
$$
f_{n} = f_{n-1}\vee q_{n}^{\prime}\leq \vee_{1}^{n-1}q_{i}\vee q_{n}
=\vee_{1}^{n}q_{i}
$$
and
$$
\tau (\vee_{1}^{n} q_{i}-f_{n})\leq ({\frac 1 2})^{k+1}+\cdot\cdot\cdot
+({\frac 1 2})^{k+n-1}<({\frac 1 2})^{k}.
$$

 For the given $\epsilon >0$ let us choose $m\in N$ such that
$\tau (I-\vee_{1}^{m}q_{i})<{\epsilon\over 2}$ and $k\in N$ such that
${\epsilon\over 2}>({\frac 1 2})^{k}$ $(>\tau (\vee_{1}^{m}g_{i}-f_{m}))$.
Then the projection $p_{\epsilon}:=f_{m}$ is that in question.

ii) Now let $e_{n}:= \wedge_{m\geq n}p_{2^{-m}}$. Then
$e^{\perp}_{n}=\vee_{m\geq n}p^{\perp}_{2^{-m}}$ and $\tau
(e_{n}^{\perp})
\leq\sum_{m\geq n}2^{-m}=2^{-n+1}$. The sequence $\{e_{n}\}$ is valid.
$\Box$

{\bf Theorem 2.} {\it Let $\cal A$ be a semifinite von Neumann algebra
containing no direct summand of type $I_{2}$ acting in the separable
Hilbert space and let
 $\mu :{\cal A}^{pr}\to [0,+\infty ]$ be the $\sigma$-finite infinite measure
and ${\cal M}_{\mu}:=\{p\in {\cal A}^{pr} :\mu (p)<+\infty\}$. Then
${\cal M}_{\mu}$ is the ideal of projections. If $\cal A$ is a finite
von Neumann algebra then the restriction
$\mu_{1}:=\mu /{\cal M}_{\mu}$ is the closed measure.}

{\it Proof.} 1). Let us prove that ${\cal M}_{\mu}$ is the ideal
of projections. It is clear that a) on ${\cal M}_{\mu}$ is fulfilled. Let
$p$, $q\in {\cal M}_{\mu}$
and $\Vert pq\Vert <1$. It is sufficient to consider the case $p$, $q$ when
$p$, $q$ are projections {\it in general position} in $H$, i.e.
$$
p\wedge q=(p\vee q-p)\wedge q = (p\vee q-q)\wedge p = 0.
\eqno{(1)}
$$
By (1), $\overline{pqH}=pH$.

i). Let us suppose first that projections $p$, $q$ are finite
with respect to $\cal A$. There exists a representation
$q=q_{1}+q_{2}+q_{3}$ (if $\cal A$ is the continuous algebra
then $q_{3}=0$) $q_{1}$, $q_{2}$,
$q_{3}\in {\cal M}$ such that the orthogonal projections
$p_{i}$ onto subspaces $\overline pq_{i}H$, $i=1$, $2$, $3$ are mutually
orthogonal and there exist the partial isometries $v_{i}\in\cal A$, $i=1$,
$2$, $3$
such that $q_{i}H$ are the initial subspaces and the final subspaces
in $(q-q_{i})H$.
The von Neumann algebra ${\cal A}^{i}$ generated by
$p_{i}$, $q_{i}$ and $v_{i}$ is direct integral of factors of type $I_{3}$.
By the construction, $p_{i}$, $q_{i}$, $v_{i}q_{i}v_{i}^{*}\in\cal M$.
By Lemma $^{(1)}$, $\mu (p_{i}\vee q_{i}-p_{i})<+\infty$ (and hence
$\mu (p_{i}\vee q_{i})<+\infty$)
if ${\cal A}^{i}$ is the type $I_{3}$ factor.
If ${\cal A}^{i}$ is the direct integral of factors, the proof
of $\mu (p_{i}\vee q_{i})<+\infty$ repeat of the proof of
Lemma $^{(1)}$. Thus
the inequality $\mu (p\vee q)=\mu (\sum_{i}p_{i}\vee q_{i})<+\infty$ is
proved.
By Lemma 5 $^{(3)}$,
$$
\mu (p\vee q)\leq (1-\Vert pq\Vert)^{-1}(\mu (p)+\mu (q)).
\eqno{(2)}
$$

ii). Let us consider now the general case of $p$, $q\in\cal M$.
Let $p_{n}\in\cal M$ be a   sequence of finite projections,
$p_{n}\nearrow p$ and let $q_{n}$ be the orthogonal projection
onto $\overline{qp_{n}H}$. The projection $q_{n}$ is finite and
$p_{n}\vee q_{n}$ $\nearrow$ $p\vee q$. The projections
$p_{n}$ and
$q_{n}$ are in the general position on the space $p_{n}\vee q_{n}H$. By (2),
$$
\mu (p\vee q)=\lim \mu (p_{n}\vee q_{n})\leq
\lim (1-\Vert p_{n}\vee q_{n}\Vert )^{-1}(\mu (p_{n})+\mu (q_{n}))\leq
$$
$$
(1-\Vert pq\Vert )^{-1})(\mu (p)+\mu (q))<+\infty.
$$
Hence $p\vee q\in {\cal M}_{\mu}$ and thus
${\cal M}_{\mu}$ is the ideal of projections.

2). Let $\{e_{n}\}$ be the sequence from Proposition 1. Then
$p\wedge e_{n}\nearrow p$, $\forall p\in {\cal A}^{pr}$. If
$\sup \{\mu_{1} (p\wedge e_{n}): n\}<+\infty$ then $p\in {\cal M}_{\mu}$.
Thus the set ${\cal M}_{\mu}$ contain any  hereditary
projection of finite $\mu$-measure. By the definition,
$\mu_{1}$ is the closed measure. $\Box$

{\bf Theorem 3.} {\it Let $\cal A$ be a finite von Neumann algebra containing
no
direct summand of type $I_{2}$ acting
in the separable Hilbert space. Then any closed measure
$\mu : {\cal M}\to [0,+\infty]$ can be extended to a $\sigma$-finite measure.}

{\it Proof.} Let $\mu :{\cal M}\to R$ be a closed measure.
By Theorem $^{(2)}$, any full finite measure can be extended
by the strong operator
topology to a unique fully finite measure on ${\cal A}^{pr}$.
Now we may assume that the measure $\mu$ is not fully finite.
Put $\mu_{1}(p):=+\infty$ for any $p\in {\cal A}^{pr}\backslash \cal M$
and  $\mu_{1}(p):=\mu (p)$,
for any $p\in\cal M$. Let us prove that the function
$\mu_{1} :{\cal A}^{pr}\to [0,+\infty ]$
is a $\sigma$-finite measure. Let $p=\sum p_{i}$ be a decomposition of
 $p\in {\cal A}^{pr}$. By a) and by the definition
of the measure, we have
 $\mu_{1} (p)=\sum_{i}\mu_{1}(p_{i})$ for any $p\in\cal M$.
 Now let $p\in {\cal A}^{pr}\backslash\cal M$. Let us assume for the moment
$\sum_{i}\mu_{1}(p_{i})<+\infty$.
By the finiteness of $\cal A$, the assumption gives us
that $p$ is the projections of finite $\mu$ measure.
 By Proposition 1, $p$ is a hereditary projection of finite
 $\mu$-measure. By the condition of the Theorem, $p\in\cal M$.
 We have the contradiction with $p\in {\cal A}^{pr}\backslash\cal M$.
Therefore $\sum_{i}\mu_{1}(p_{i})=+\infty =\mu_{1}(p)$.
Let $\{e_{n}\}$ be  the sequence from Proposition 1.
Hence $\mu_{1}$ is a $\sigma$-finite measure. $\Box$

\bigskip
\noindent
{\bf ACKNOWLEDGMENT}.

The research supported
by the grant Min. Obrazovaniya Rossii E00-1.0-172

\bigskip
\noindent
{\bf REFERENCES}

1. Lugovaja G.D. and Sherstnev A.N. (1980). On the Gleasons theorem for
unbounded
measures. {\it Izvestija VUZov. Matematika}. n.12, 30-32. [in Russian].

2. Matvejchuk M.S. (1981). Description of finite measures on semifinite
algebras.
    {\it Functional Anal. i Prilozhen.}, vol.15, n.3, 41-53. [in Russian],
    Engl. transl.: Funct.Anal.Appl. 1981, vol.15, n.3, 187-197. MR\#
84h:46088.

3. Matvejchuk M.S. (1981). Finite signed measures on von Neumann algebras.
    {\it Constructive theory of functions and functional analysis}, Kazan Cos.
Univ. Kazan. vol.III, 55-63. [in Russian], MR\# 83i:46075




\end{document}